\def\AL{\hbox{\rm A$_L$\hskip 2pt}}
\def\grad{\overline\nabla_X f}

\def\IC{{\bf C}} 
\def\TheMagstep{\magstep1}	
\def\PaperSize{letter}		

\magnification=\magstep1

\let\:=\colon  \let\To=\longrightarrow
   \let\?=\overline 
 \let\To=\longrightarrow

\let\Sum=\sum \def\sum{\Sum\nolimits}

\def\IC{{\bf C}} 
\def\IP{{\bf P}}

\def\pd #1#2{\partial#1/\partial#2}

\def\and{\hbox{ and }}
\def\Wf{\hbox{\rm W$_f$\hskip 2pt}}
	
\def\Af{\hbox{\rm A$_f$\hskip 2pt}}

\def\AF{\hbox{\rm A$_F$\hskip 2pt}}

\def\DONE{*!*}
\def\NextDef #1 {\def\NextOne{#1}%
 \ifx\NextOne\DONE\let\next\relax
 \else\expandafter\xdef\csname#1\endcsname{\TheOp}
  \let\next\NextDef
 \fi \next}
\def\TheOp{\mathop{\rm\NextOne}}
 \NextDef 
  Projan Supp Proj Sym Spec Hom cod Ker dist
 *!*
\def\TheOp{{\cal\NextOne}}
\NextDef 
  E F G H I J M N O R S
 *!*
\def\TheOp{\hbox{\rm\NextOne}}
\NextDef 
 A ICIS 
 *!*

\def\item#1 {\par\indent\indent\indent\indent \hangindent4\parindent
 \llap{\rm (#1)\enspace}\ignorespaces}
 \def\inpart#1 {{\rm (#1)\enspace}\ignorespaces}
 \def\part {\par\inpart}

\catcode`\@=11		

\def\vfootnote#1{\insert\footins\bgroup
 \eightpoint 
 \interlinepenalty\interfootnotelinepenalty
  \splittopskip\ht\strutbox 
  \splitmaxdepth\dp\strutbox \floatingpenalty\@MM
  \leftskip\z@skip \rightskip\z@skip \spaceskip\z@skip \xspaceskip\z@skip
  \textindent{#1}\footstrut\futurelet\next\fo@t}

\def\p.{p.\penalty\@M \thinspace}
\def\pp.{pp.\penalty\@M \thinspace}
\def\(#1){{\rm(#1)}}\let\leftp=(
\def\activeleftp{\catcode`\(=\active}
{\activeleftp\gdef({\ifmmode\let\next=\leftp \else\let\next=\(\fi\next}}

\def\sct#1\par
  {\removelastskip\vskip0pt plus2\normalbaselineskip \penalty-250 
  \vskip0pt plus-2\normalbaselineskip \bigskip
  \centerline{\smc #1}\medskip}

\newcount\sctno \sctno=0
\def\sctn{\advance\sctno by 1 
 \sct\number\sctno.\quad\ignorespaces}

\def\dno#1${\eqno\hbox{\rm(\number\sctno.#1)}$}
\def\Cs#1){\unskip~{\rm(\number\sctno.#1)}}

\def\proclaim#1 #2 {\medbreak
  {\bf#1 (\number\sctno.#2)}\enspace\bgroup\activeleftp
\it}
\def\endproclaim{\par\egroup\medskip}
\def\pf{\endproclaim{\bf Proof.}\enspace}
\def\lem{\proclaim Lemma } \def\prp{\proclaim Proposition }
\def\cor{\proclaim Corollary }	\def\thm{\proclaim Theorem }
\def\rmk#1 {\medbreak {\bf Remark (\number\sctno.#1)}\enspace}
\def\eg#1 {\medbreak {\bf Example (\number\sctno.#1)}\enspace}

\parskip=0pt plus 1.75pt \parindent10pt
\hsize29pc
\vsize44pc
\abovedisplayskip6pt plus6pt minus2pt
\belowdisplayskip6pt plus6pt minus3pt

\def\TRUE{TRUE}	
\ifx\DoublepageOutput\TRUE \def\TheMagstep{\magstep0} \fi
\mag=\TheMagstep

\newskip\vadjustskip \vadjustskip=0.5\normalbaselineskip
\def\centertext
 {\hoffset=\pgwidth \advance\hoffset-\hsize
  \advance\hoffset-2truein \divide\hoffset by 2\relax
  \voffset=\pgheight \advance\voffset-\vsize
  \advance\voffset-2truein \divide\voffset by 2\relax
  \advance\voffset\vadjustskip
 }
\newdimen\pgwidth\newdimen\pgheight
\def\letter{letter}\def\AFour{AFour}
\ifx\PaperSize\letter
 \pgwidth=8.5truein \pgheight=11truein 
 \message{- Got a paper size of letter.  }\centertext 
\fi
\ifx\PaperSize\AFour
 \pgwidth=210truemm \pgheight=297truemm 
 \message{- Got a paper size of AFour.  }\centertext
\fi

 \newdimen\fullhsize \newbox\leftcolumn
 \def\fulline{\hbox to \fullhsize}
\def\doublepageoutput
{\let\lr=L
 \output={\if L\lr
          \global\setbox\leftcolumn=\columnbox \global\let\lr=R%
        \else \doubleformat \global\let\lr=L\fi
        \ifnum\outputpenalty>-20000 \else\dosupereject\fi}%
 \def\doubleformat{\shipout\vbox{%
        \fulline{\hfil\hfil\box\leftcolumn\hfil\columnbox\hfil\hfil}%
				}%
		  }%
 \def\columnbox{\vbox
   {\makeheadline\pagebody\makefootline\advancepageno}%
   }%
 \fullhsize=\pgheight \hoffset=-1truein
 \voffset=\pgwidth \advance\voffset-\vsize
  \advance\voffset-2truein \divide\voffset by 2
  \advance\voffset\vadjustskip
 \let\firstheadline=\hfil
 
}
\ifx\DoublepageOutput\TRUE \doublepageoutput \fi

 \font\twelvebf=cmbx12		
 \font\smc=cmcsc10		

\def\eightpoint{\eightpointfonts
 \setbox\strutbox\hbox{\vrule height7\p@ depth2\p@ width\z@}%
 \eightpointparameters\eightpointfamilies
 \normalbaselines\rm
 }
\def\eightpointparameters{%
 \normalbaselineskip9\p@
 \abovedisplayskip9\p@ plus2.4\p@ minus6.2\p@
 \belowdisplayskip9\p@ plus2.4\p@ minus6.2\p@
 \abovedisplayshortskip\z@ plus2.4\p@
 \belowdisplayshortskip5.6\p@ plus2.4\p@ minus3.2\p@
 }
\newfam\smcfam
\def\eightpointfonts{%
 \font\eightrm=cmr8 \font\sixrm=cmr6
 \font\eightbf=cmbx8 \font\sixbf=cmbx6
 \font\eightit=cmti8 
 \font\eightsmc=cmcsc8
 \font\eighti=cmmi8 \font\sixi=cmmi6
 \font\eightsy=cmsy8 \font\sixsy=cmsy6
 \font\eightsl=cmsl8 \font\eighttt=cmtt8}
\def\eightpointfamilies{%
 \textfont\z@\eightrm \scriptfont\z@\sixrm  \scriptscriptfont\z@\fiverm
 \textfont\@ne\eighti \scriptfont\@ne\sixi  \scriptscriptfont\@ne\fivei
 \textfont\tw@\eightsy \scriptfont\tw@\sixsy \scriptscriptfont\tw@\fivesy
 \textfont\thr@@\tenex \scriptfont\thr@@\tenex\scriptscriptfont\thr@@\tenex
 \textfont\itfam\eightit	\def\it{\fam\itfam\eightit}%
 \textfont\slfam\eightsl	\def\sl{\fam\slfam\eightsl}%
 \textfont\ttfam\eighttt	\def\tt{\fam\ttfam\eighttt}%
 \textfont\smcfam\eightsmc	\def\smc{\fam\smcfam\eightsmc}%
 \textfont\bffam\eightbf \scriptfont\bffam\sixbf
   \scriptscriptfont\bffam\fivebf	\def\bf{\fam\bffam\eightbf}%
 \def\rm{\fam0\eightrm}%
 }

\def\today{\ifcase\month\or	
 January\or February\or March\or April\or May\or June\or
 July\or August\or September\or October\or November\or December\fi
 \space\number\day, \number\year}
\nopagenumbers
\headline={%
  \ifnum\pageno=1\firstheadline
  \else
    \ifodd\pageno\oddheadline
    \else\evenheadline\fi
  \fi
}
\let\firstheadline\hfill
\def\oddheadline{
 \hfil\headtitle\hfil\llap{\folio}}
\def\evenheadline{\eightpoint\rlap{\folio}
 \hfil\author\hfil\llap{\today}}
\def\headtitle{\title}

 \newcount\refno \refno=0	 \def\NoKey{*!*}
 \def\MakeKey{\advance\refno by 1 \expandafter\xdef
  \csname\TheKey\endcsname{{\number\refno}}\NextKey}
 \def\NextKey#1 {\def\TheKey{#1}\ifx\TheKey\NoKey\let\next\relax
  \else\let\next\MakeKey \fi \next}
 \def\RefKeys #1\endRefKeys{\expandafter\NextKey #1 *!* }
\def\SetRef#1 #2,#3\par{%
 \hang\llap{[\csname#1\endcsname]\enspace}%
  \ignorespaces{\smc #2,}
  \ignorespaces#3\unskip.\endgraf
 }
 \newbox\keybox \setbox\keybox=\hbox{[8]\enspace}
 \newdimen\keyindent \keyindent=\wd\keybox
\def\references{
  \bgroup   \frenchspacing   \eightpoint
   \parindent=\keyindent  \parskip=\smallskipamount
   \everypar={\SetRef}}
\def\endreferences{\egroup}

 \def\serial#1#2{\expandafter\def\csname#1\endcsname ##1 ##2 ##3
  {\unskip\ #2 {\bf##1} (##2), ##3}}
 \serial{ajm}{Amer. J. Math.}
  \serial {aif} {Ann. Inst. Fourier}
 \serial{asens}{Ann. Scient. \'Ec. Norm. Sup.}
 \serial{comp}{Compositio Math.}
 \serial{conm}{Contemp. Math.}
 \serial{crasp}{C. R. Acad. Sci. Paris}
 \serial{dlnpam}{Dekker Lecture Notes in Pure and Applied Math.}
 \serial{faa}{Funct. Anal. Appl.}
 \serial{invent}{Invent. Math.}
 \serial{ma}{Math. Ann.}
 \serial{mpcps}{Math. Proc. Camb. Phil. Soc.}
 \serial{ja}{J. Algebra}
 \serial{splm}{Springer Lecture Notes in Math.}
 \serial{tams}{Trans. Amer. Math. Soc.}

\def\UThin{\penalty\@M \thinspace\ignorespaces}
\def\relaxnext@{\let\next\relax}
\def\cite#1{\relaxnext@
 \def\nextiii@##1,##2\end@{\unskip\space{\rm[\SetKey{##1},\let~=\UThin##2]}}%
 \in@,{#1}\ifin@\def\next{\nextiii@#1\end@}\else
 \def\next{{\rm[\SetKey{#1}]}}\fi\next}
\newif\ifin@
\def\in@#1#2{\def\in@@##1#1##2##3\in@@
 {\ifx\in@##2\in@false\else\in@true\fi}%
 \in@@#2#1\in@\in@@}
\def\SetKey#1{{\bf\csname#1\endcsname}}

\catcode`\@=12 
\def\title{ The Multiplicity Polar Theorem and Isolated Singularities}
\def\author{Terence Gaffney}
\RefKeys  Br BLS BMPS  B-R E-GZ1 E-GZ2 E-GZ3 EGZS
 F FM  G-1 G-2 G-3 G-4 G-5 G-6 G-7 G-8  GK   G
K-T KT1 L Le1 Le2   MP  
Ma Ma1 Ma-2  T-2 
 \endRefKeys

\def\topstuff{\leavevmode
 \bigskip\bigskip
 \centerline{\twelvebf \title}
 \bigskip
 \centerline{\author}
 \medskip\centerline{\today}
\bigskip\bigskip}
\topstuff
\sct Introduction

Recently there has been much work done in the study of isolated singularities. 
The singularities may be those of sets, functions, vector fields or differential
forms. Examples can be found in the references at the end of this paper. In these works the authors have developed algebraic invariants which describe the
isolated singularity. In this paper, we show how many of these invariants can be computed by relating them to the multiplicity of a pair of modules, by means of the
multiplicity polar theorem. The method of application is the same in all cases. A deformation is constructed in which the isolated singularity breaks into
the simplest possible singularities, which are well understood, and the multiplicity polar theorem relates the multiplicity of the pair of modules in the
original situation to the multiplicities at the simple points.

Our answer involves a pair of modules because the the multiplicity is only defined for modules and ideals of finite colength, and in many of the situations we want to describe the modules do not have finite colength.

The solution to this problem is to use pairs of modules, $(M,N)$, where $M\subset N\subset {\cal O}^p_{X,x}$; here $M$ has finite colength in $N$, $N$ is the "least complicated" module in ${\cal O}^p_{X,x}$  consistent with the geometry of $X$ that makes the colength of $M$ infinite. The difference between $M$ and $N$, which is measured by the multiplicity of the pair, should reflect only the contribution of the geometry at the point $x$. 

This geometry can be related to other invariants of $X$ by using an $N$ which reflects some of the geometry of $X$ at the point, then using the additivity property of multiplicity to relate the multiplicities of the pairs. This is illustrated in propositions 2.5 and 2.6.

The first application is to the study of the Euler obstruction.

 The Euler Obstruction is the basic tool used by MacPherson in constructing a Chern class on singular spaces. (Cf.\cite{FM}, \cite {Br} for
details.)

Brasselet,  L\^e,  Seade, \cite{BLS} found a formula for the Euler obstruction using a general linear form $L$. Their formula is:
 $${\rm Eu}_X(0)=\sum_i  \chi(V_i\cap  B_\varepsilon \cap L^{-1}(t_0)) \cdot {\rm Eu}_X(V_i),$$
where $\ B_\varepsilon$ is a small ball around $0$ in $C^N$, $t_0  
\in
C
\setminus \{ 0 \}$ is sufficiently near $\{ 0 \}$ and $Eu_X(V_i)$ is the
Euler
obstruction of $X$ at any point of the stratum $V_i$, $\{V_i\}$ a Whitney stratification.

Then Brasselet, Massey, Parameswaran and  Seade (\cite{BMPS})   asked ``Suppose a function $f$ defined on $X$ is used in the formula instead
of a general linear form;  how far off is the number obtained from the Euler obstruction?" This difference is called the defect of $f$ and
is denoted $D_{f, X}(0)$.

One reason that the question is interesting, is that the answer gives a measure of the singularity of $f$ 
which is independent of the singularities of $X$--a generic linear function 
is the ``least singular" function on $X$. So the singularities of $f$ beyond the singularities of $L$, are those which are due to $f$, 
not $X$. In the case where $f$ has an isolated singularity, the defect also has an interpretation as the obstruction for constructing a section of
the Nash bundle over the Nash transform of $X$ that extends a lifting of $\grad$,
the complex conjugate of the gradient vector field of $f$ on $X$ (\cite{BMPS}).

Suppose $f$ has an isolated singularity on $X$, and $T^*(X)$ is the closure of the conormal vectors 
to the smooth part of $X$. Denote the image of $df$ over $X$  by $\Gamma (df)$. Then BMPS showed that the defect, up to a sign, 
is the intersection multiplicity $\Gamma (df)\cdot T^*(X)$, and they computed this in terms of the 
L\^e-Vogel numbers of a certain sheaf \cite{Ma-2}.

The theory of the multiplicity of pairs of modules provides another approach to this computation.
 In this paper we will 
 describe, in terms of the theory of integral closure of modules, what it means for $f$ to have an isolated singularity on $X$. 
This description will naturally 
associate a pair of modules to $f$; the multiplicity of the pair will be the desired intersection number. Led by evidence from the ICIS case, we will associate a related pair of modules to $f$ and  to a generic linear function $L$. 
We will then show that the desired intersection multiplicity is the
 difference of the multiplicities
of the associated pairs of modules. Our constructions also apply to a 1-form; then the difference can be interpreted as the index of the 1
form. The same computation also is relevant to the theory of D-modules. Using work of L\^e, we can get a formula for the complex
dimension of the space of vanishing cycles of a sheaf of D-modules $M$ relative to a function $f$ at $0$, in terms of the multiplicities of the 
characteristic cycle of $M$ and the multiplicities of a pair of modules associated to $f$ and a generic linear form.

The same computation also applies to work of Massey, where we use it to get a formula for the relative cohomology  of the Milnor fiber of $f$ where $f$
has an isolated singularity on a complex analytic set with possibly non-isolated singularities. In the special case in which $X$ is a hypersurface, with
singular locus an ICIS, and Morse transverse singularity, we are able to eliminate the dependence on the linear form from Massey's result.

In our final application we return to the study of the \Af condition on a family of spaces with isolated singularities, and there we use the information
from a generic linear function on the family to refine the invariants from \cite{G-6}, removing some of the dependence on the family. In the event that, in addition, the family satisfies a continuity condition, we close the circle of ideas in the paper by showing that the \Af condition holds if and only if the defect of the function at each member of the family is constant.

The main technical tool, as mentioned earlier, is the Multiplicity-Polar Theorem, which is discussed in the first section. To show how this theorem is used
in the following sections, we give a new proof in the context of analytic spaces of a result of Buchsbauum and Rim relating the colength of a module and its
ideal of maximal minors.

It is a pleasure to thank David Massey for introducing me to \cite{BMPS}, and for helpful conversations in developing this paper, particularly the
material around Theorem 3.3. I also thank the organizers of the Winter School in singularities at Luminy. Their invitation was a spur to completing the work.

Since the multiplicity-polar theorem plays such an important role in this paper, I also thank Steven Kleiman for conversations that led up to its proof.

\sctn {Background on the theory of integral closure of modules}

In this paper we work with complex analytic sets and maps. Let $\O_X$ 
denote the structure sheaf on a complex analytic space $X$.

The study of what it means for a function to have an isolated singularity on a singular space in the next section,
depends on the behavior of
limiting  tangent hyperplanes. The key tool for studying these limits is the theory of integral closure of modules, which we now introduce.

Suppose $(X,x)$ is a complex
analytic germ, $M$ a submodule of
$\O_{X,x}^{p}$. Then $h\in \O_{X,x}^{p}$ is in the {\it{integral closure of
$M$}}, denoted $\?{M}$,
if and only if for all
$\phi :\left({\IC},0\right) \rightarrow \left(X,x\right)$, $ h \circ \phi
\in (\phi^{*}M){\O}_{1}$.

If $N$ is a submodule of $M$ and $\?M=\?N$ we say that $N$ is a {\it reduction} of $M$.
(For more details cf. \cite {G-2}).

If a module $M$ has finite colength in $\O^p_{X,x}$, it is possible to 
attach a number 
to the module, its Buchsbaum-Rim multiplicity (\cite{B-R}). We can also 
define the multiplicity
of a pair of modules $M\subset N$, $M$ of finite colength in $N$, as well, 
even if $N$ does not have finite colength in $\O^p_X$. 

In studying the geometry of singular spaces, it is natural to study pairs of modules. In dealing with non-isolated singularities, the
modules that describe the geometry have non-finite colength, so their multiplicity is not defined. Instead it is possible to define a
deceasing sequence of modules, each with finite colength inside its predecessor, when restricted to a suitable complentary plane. Each pair
controls the geometry in a particular codimension. For an example of how to construct this sequence see \cite{G-5}.

It is also interesting to consider families with specified singularities, for example familes of hypersurfaces  where the singular set 
is a family of ICIS singularities defined by $I$. Here the members of the family
 are the pairs $(J(f_t), I(t))$ where $t$ is the parameter
of the family. By specifying the pair, we capture the restriction we have in mind on the possible deformations of $f$. Since the possible
deformations become finite again, the multiplicity is well defined. For
examples of this cf. \cite {G-6} and \cite {G-7}.

We recall how to construct the multiplicity following the approach
 of Kleiman and Thorup (\cite{K-T}).

Given a submodule $M$ of a free $\O_X$ module $F$ of rank $p$, we can 
associate a subalgebra ${\cal R }(M)$ of the symmetric
$\O_X$ algebra on $p$ generators. This is known as the Rees algebra of 
$M$. If $(m_1,\dots,m_p)$ is an element of $M$
 then $\Sum m_iT_i$ is the corresponding element of 
${\cal R}(M)$. Then $\Projan ({\cal R} (M))$, the projective analytic 
spectrum of ${\cal R}(M)$ is the closure of 
the projectivised row spaces of $M$ at points where the rank
of a matrix of generators of $M$ is maximal. Denote the projection to $X$ 
by $c$, or by $c_M$ where there is ambiguity.

 If $M$ is a submodule of $N$ or $h$ is a section of $N$, then
 $h$ and  $M$  generate ideals on $\Projan{\cal R}(N)$; denote them
by $\rho(h)$ and $\rho({\cal M})$. If we can express $h$ in terms of a set 
of generators $\{n_i\}$ of $N$ as $\Sum g_in_i$, 
then in the chart in which $T_1\ne 0$, we can express a generator of 
$\rho(h)$ by $\Sum g_iT_i/T_1$. 
Having defined the ideal sheaf
$\rho({\cal M})$, we blow up by it.

On the blowup $B_{\rho({\cal M})}(\Projan {\cal R}(N))$ we have two 
tautological bundles, one the pullback of 
the bundle on $\Projan{\cal R}(N)$,
 the other coming from $\Projan{\cal R}(M)$; denote the corresponding 
Chern classes by $l_M$ and $l_N$, and denote the exceptional divisor by
$ D_{M,N}$. 
Suppose the generic 
rank of $N$ (and hence of $M$) is $e$.
 Then 
the multiplicity of a pair of modules $M,N$ is:

$$e(M,N)=\Sum_{j=0}^{d+e-2} \int D_{M,N} \cdot l^{d+e-2-j}_{M}\cdot 
l^j_{N}.$$

Kleiman and Thorup show that this multiplicity is well defined at $x\in X$ as long as 
$\bar M=\bar N$ on a deleted neighborhood of $x$. This condition implies that $D_{M,N}$ lies in the fiber over $x$ hence is compact.

If $M\subset \O^p_X$ is of finite colength, then the multiplicity of $M$ 
is the multiplicity of the pair $(M,\O^p_X)$.

If $\O_{X^d,x}$ is Cohen-Macauley, and $M$ has $d+p-1$ generators then 
there is a useful relation between $M$ and its ideal of maximal minors; 
the multiplicity of $M$ is the colength of $M$, is the
 colength of the ideal of maximal minors,
 by some theorems of Buchsbaum and Rim \cite{B-R}, 2.4 \p.207, 4.3 and
4.5 \p.223. As an illustration of the power of the multiplicity polar theorem 
we will give a proof of this theorem in the context of
analytic geometry.

We next develop the notion of polar varieties which is the other term in 
the multiplicity polar theorem.

Assume we have a module $M$ which is a submodule of a free module on $X^d$ 
an equidimensional,  analytic space, reduced off a nowhere dense subset of 
$X$,  and assume the generic rank of $M$ is $e$ on each 
component of $X$. The hypothesis on the equidimensionality of $X$ and on 
the rank of $M$ ensures that $\Projan{\cal R}(M)$ 
is equidimensional of 
dimension $d+e-1$. Note that $\Projan{\cal R}(M)$ can be embedded in 
$X\times \IP^{r-1}$, 
provided we can chose a set of generators of $M$ 
with $r$ elements.
The {\it polar variety of codimension $k$} of $M$ in $X$ denoted 
$\Gamma_k(M)$ is constructed by intersecting $\Projan{\cal R}(M)$
 with $X\times H_{e+k-1}$ where 
$H_{e+k-1}$ is a general plane of codimension $e+k-1$, then projecting to 
$X$. 
When we consider $M$ as part of a pair of modules $M,N$, where the generic 
rank of $M$ is the same as the generic rank of $N$,
 then other polar varieties become interesting as well. In brief, we can 
intersect
 $B_{\rho({\cal M})}(\Projan {\cal R}(N))\subset X\times 
\IP^{N-1}\times\IP^{p-1}$ with a mixture of hyperplanes from the two 
projective spaces
which are factors of the space in which the blowup is embedded. We can 
then push these intersections down to 
$\Projan {\cal R}(N)$ or
$X$ as is convenient, getting mixed polar varieties in $\Projan {\cal 
R}(N)$ or in $X$. These mixed varieties play an 
important role in the proof of the multiplicity-polar theorem, the theorem 
we next describe.

Setup: We suppose we have families  of modules $M\subset  N$, $M$ and $N$
 submodules of a free module $F$ of rank $p$
 on an equidimensional family of spaces with equidimensional
 fibers ${\cal X}^{d+k}$, ${\cal X}$ a family over a smooth base
$Y^k$. We assume that the generic rank of $M$, $N$ is $e\le p$.  Let 
$P(M)$ denote $\Projan {\cal R}(M)$, $\pi_M$
 the projection to ${\cal X}$.  
let $C(M)$ denote the locus of points where $M$ is not free, ie. the 
points where the rank of $M$ is less 
than $e$, $C(\Projan {\cal R}(M))$
its inverse image under $\pi_M$, $C({\cal M})$ the cosupport of 
$\rho({\cal M})$ in $P(\Projan {\cal R}(N))$.

We will be interested in computing the change in the multiplicity of the 
pair $(M,N)$, denoted  $\Delta (e(M,N))$. We will
assume that the integral closures of $M$ and $N$ agree  off a set $C$ of 
dimension $k$ which is finite over $Y$, and assume
we are working on a sufficently small neighborhood of the origin, that 
every component of $C$ contains the origin in its
closure. Then $e(M,N,y)$ is the sum of the multiplicities of the pair at 
all points in the fiber of $C$ over $y$, and $\Delta (e(M,N))$
is the change in this number from $0$ to a generic value of $y$. If we 
have a set $S$ which is finite over $Y$, then we can project $S$ to $Y$, 
and the degree
of the branched cover at $0$ is ${\rm  mult}_y S$. (Of course, this is 
just the number of points in the fiber of $S$ over our generic $y$.)

We can now state our theorem.

\thm 1 Suppose in the above setup we have that $\?M=\?N$ off a set $C$ of 
dimension $k$ which is finite over $Y$. 
Suppose further
that $C(\Projan {\cal R}(M))(0)=C(\Projan {\cal R}(M(0)))$ except possiby 
at the points which project to $0\in {\cal X}(0)$. 
Then, for $y$ a 
generic point of $Y$,

$$\Delta (e(M,N))={\rm  mult}_y \Gamma_{d}(M)-{\rm mult}_y 
\Gamma_{d}(N).$$

\pf The proof in the ideal case appears in \cite{G-6}; the general proof 
will appear in \cite{G-8}. 

The first condition in the theorem ensures the multiplicity is defined fiberwise, while the second implies that the conormal geometry of $M$ in the special fiber away from $0$ is not different from that of the general fiber.

The following application appears in \cite {G-7} and will be used to prove the connection between 
the colenth of a module and its associated ideal of maximal minors.

\thm 2 Given $M$ a submodule of  $\O^p_{X,x}$, $X^d$ equidimensional, choose
$d+p-1$ elements which generate a reduction $K$ of $M$. Denote the matrix 
whose columns are the $d+p-1$ elements by $[K]$; $[K]$ induces
a section of $\Hom \hskip 2pt(\IC^{d+p-1},\IC^p)$ which is a trivial 
bundle over $X$. Stratify $\Hom \hskip 2pt(\IC^{d+p-1},\IC^p)$ by rank. 
Let $[\epsilon]$ denote a $p\times (d+p-1)$ matrix, whose entries are 
small, generic constants. Then, on a suitable neighborhood $U$ of $x$ 
the section of 
$\Hom \hskip 2pt(\IC^{d+p-1},\IC^p)$ induced from $[K]+[\epsilon]$ has at 
most kernel rank 1, is transverse to the rank stratification, 
and the number of points where the kernel rank is 1 is $e(M)$.
\endproclaim

 We now use the underlying idea of this proof to make the connection between the colength of
 the module and the ideal of its maximal minors.

\thm 3 Suppose $M$ a submodule of finite colength of $\O^p_{X^d,x}$ with $d+p-1$ generators, $\O^p_{X^d,x}$ Cohen-Macaulay, and $J(M)$ the
ideal of maximal minors of $M$. Then the colength of $M$ and the colength of $J(M)$ are the same.

\pf Form a family of modules ${\cal M}$ on $X\times \IC$ with a matrix of generators $[M]+t[\epsilon ]$, 
where $[M]$ is a matrix of generators for $M$ with
$d+p-1$ columns and $[\epsilon]$ the matrix of generic constants. The proof of Theorem 1.2 uses the multiplicity polar theorem
to show that the colength of $M$ is the number of points where for generic t, the rank of $[M]+t[\epsilon ]$ is less than maximal. 

Now consider the analytic set defined by $J({\cal M})$ on  $X\times \IC$. The fiber over a generic point $t$ close to zero is again just the
points where the rank of $[M]+t[\epsilon ]$ is less than maximal. Since  $\O^p_{X^d,x}$ is Cohen-Macaulay, so is $V({\cal J}( M))$, hence
the degree of the projection to the base $\IC$ at $(x,0)$ is just the colength of $J(M)$. So the colength of $J(M)=$ the number of points
in a generic fiber $=$ the colength of $M$.

\sctn {Isolated singularities and the defect.}

As mentioned in the introduction, our notion of isolated singularity depends on  the behavior of limiting 
tangent hyperplanes.

    It
is natural to study these limits by making them part of a space; the relative conormal space of $f$ denoted $C(f)$, is defined
 by considering 
the set of tangent hyperplanes to the fibers of $f$ at points where $X$ and $f$ are smooth, and taking the closure of this
 set in $X\times\IP^{n-1}$.  (The conormal space of $X$, $C(X)$, is the relative conormal space of the zero map, viewed as a 
map from $X$ to 
$0$.)

If $X$ is defined by $F:\IC^N\to\IC^p$, then the Jacobian module of $X$ denoted $JM(X)$, is the submodule of $\O^p_X$ generated by the
partial derivatives of 
$F$. Given a function $f$ defined on $X$,  form the $p+1$ by $n$ matrix $D(F,f)$
by augmenting the Jacobian matrix $DF$ at the bottom with the gradient
$df$. Call the submodule of the free module $\O_X^{p+1}$, generated by
the columns of $D(F,f)$, the {\it augmented Jacobian module\/} and
denote it by $JM(X,f)$.
At a point $x$ where $X$ and $f$ are smooth, the tangent hyperplanes to the fiber of $f$ at $x$ can be identified
with the projectivised row space of the augmented Jacobian matrix, so $C(f)=\Projan(R(JM(X,f)))$, while $C(X)=\Projan(R(JM(X)))$.

There are different possible definitions for the singular points of $f$ (cf. \cite{Ma-2}  for a discussion of the possibilties.)
For this section we say that $f$ has an isolated singularity at the origin if $\Gamma (df)\cdot T^*(X)$ is an isolated point.

The plan then of this section is as follows. We first see what it means for a linear function to have no singularity at the origin, then for
the same to be true of any function, in integral closure terms. We then describe what it means in integral closure terms for $f$ to have an
isolated singularity. This description will naturally involve a pair of modules. We will then construct a suitable deformation of $f$,
which when coupled with the multiplicity polar theorem, will compute the intersection multiplicity of $\Gamma (df)\cdot T^*(X)$. Up to a
sign, this number is the defect.

To state a first result we need some more notation.

Given an analytic map germ $g\:(\IC^n,0)\to (\IC^l,0)$, let
$JM(X)_g$ denote the submodule of $JM(X)$ generated by the ``partials''
$\pd Fv$ for all vector fields $v$ on $\IC^n$ tangent to the fibers of
$g$,  call
$JM(X)_g$ the {\it relative Jacobian module\/} with respect to $g$.
For example, if $g$ is the projection onto the space of the last $l$
variables of $\IC^n$, then $JM(X)_g$ is simply the submodule generated
by all the partial derivatives of $F$ with respect to the first $n-l$
variables.

The following lemma describes the limit tangent hyperplanes in terms of Jacobian
modules.  

\lem1  A hyperplane $H$, defined by the vanishing of a linear
function $h\:\IC^n\to \IC$, is a limit tangent hyperplane of $(X,0)$ if
and only if $JM(X)_h$ is not a reduction of $JM(X)$.

\pf Cf. Theorem~2.4 of \cite{G-1} or lemma 4.1 of \cite{GK}.

It will be  helpful to reformulate this condition in a way that separates more clearly the infinitesimal data from $X$ and $h$.

\lem2 In the situation of Lemma 2.1,  a hyperplane $H$, defined by the vanishing of a linear
function $h\:\IC^n\to \IC$, is a limit tangent hyperplane of $(X,0)$ if
and only if $JM(X,h)$ is not a reduction of $JM(X)\oplus \O_X$.

\pf It suffices to show that  $JM(X,h)$ is  a reduction of $JM(X)\oplus \O_X$ if and only if  $JM(X)_h$ is a reduction of $JM(X)$.
Suppose $JM(X,h)$ is  a reduction of $JM(X)\oplus \O_X$, then  $\?{JM(X,h)}$ contains $JM(X)\oplus 0$. Restricting to curves, this implies
$\?{JM(X)_h}$ contains $JM(X)$.

Suppose  $JM(X)_h$ is a reduction of $JM(X)$. Then $\?{JM(X,h)}$ contains
 
\noindent $JM(X)\oplus 0$. Let $v$ be a vector so that $D(h)(v)\ne0$
Then $\?{JM(X,h)}$ contains $D(F,h)(v)$ and $JM(X)\oplus 0$, so it contains $JM(X)\oplus \O_X$.

Now we are ready to say what it means for a  function to be non-singular. Since we can think of $df$ as a family of hyperplanes which
miss $T^*(X)$, it is not surprising to get a similar criterion.

\lem3 A function $f$ has  no singularity at $0$ (ie. $df$ misses $T^*{(X)}$ ) if and only 
if  $JM(X,f)$ is  a reduction of $JM(X)\oplus \O_X$.

\pf Let $h$ be a linear function such that $h=df(0)$. Then $f$ has no singularity at $0$ iff  $JM(X,h)$ is  a reduction of $JM(X)\oplus
\O_X$. So we need to show that  $JM(X,h)$ is  a reduction of $JM(X)\oplus
\O_X$ if and only if $JM(X,f)$ is  a reduction of $JM(X)\oplus \O_X$.

Suppose  $JM(X,h)$  is a reduction of $JM(X)\oplus
\O_X$. Then since $m(JM(X)\oplus \O_X)+JM(X,f)$ contains $JM(X,h)$ it follows that
$$\?{m(JM(X)\oplus \O_X)+JM(X,f)}\supset JM(X)\oplus \O_X.$$
Then by the integral closure form of Nakayama's lemma,  $\?{JM(X,f)}\supset JM(X)\oplus \O_X$.

The converse is similar, since  $m(JM(X)\oplus \O_X)+JM(X,h)$ contains $JM(X,f)$.

If $f$ has an isolated singularity, then we expect $JM(X,f)$ to be a reduction of $JM(X)\oplus \O_X$ except at the origin. This suggests
using the pair of modules $JM(X,f), JM(X)\oplus \O_X$, as the integral closures will be the same on a deleted neighborhood of $x$, so the multiplicity will be well defined.  As we shall see, from looking at the ICIS case, it is helpful in understanding this multiplicity to use the pair $({JM(X,f)}, H_{d-1}(X)\oplus \O_X)$, where
$H_{d-1}(X)$ by definition consists of all elements of $\O^{p}_X$ which are in $\?{JM(X)}$ except at the origin. In some cases, it is easier to see what $H_{d-1}(X)$ is than $JM(X)$, and the theory with
$H_{d-1}(X)$ makes the role of a generic linear function clearer. First we give our characterization of isolated singularities 
in integral closure terms, then our next
three results explicate these points.

\prp 4  A function germ $f$ has an isolated singularity at $x\in X$ if and only if  $JM(X,f)$ is a reduction of $ H_{d-1}(X)\oplus
\O_X$ except possibly at $x$.

\pf Let $U$ be a neighborhood of $x$ in $X$ such that either $df$ misses $T^*{(X)}$ on $U-x$
 or  $JM(X,f)$ is a reduction of $ H_{d-1}(X)\oplus
\O_X$ except possibly at $x$ and $ H_{d-1}(X)=\?{JM(X)}$ except possibly at $x$. Let $x'\in U$, $x'\ne x$, then in the first case, $f$ has
no singularity at
$x'$, so  by 2.3 $JM(X,f)$ is  a reduction of $JM(X)\oplus \O_X$ at $x'$. Since by choice of $U$, $H_{d-1}(X)=\?{JM(X)}$, $JM(X,f)$ is a
reduction of $ H_{d-1}(X)\oplus
\O_X$. In the second case  $JM(X,f)$ is a reduction of $ H_{d-1}(X)\oplus
\O_X$ implies  $JM(X,f)$ is  a reduction of $JM(X)\oplus \O_X$, hence by 2.3, $f$ has no singularity at $x'$.

\prp 5 In the above setup, 
$$e(JM(X,f), JM(X)\oplus \O_X,x)$$
$$=e(JM(X,f), H_{d-1}(X)\oplus \O_X,x)-e(JM(X,L), H_{d-1}(X)\oplus \O_X,x)$$

\noindent where $L$ is a linear function with no singularity at $x$.

\pf The proof is based on a fundamental result due to Kleiman and Thorup--the principle of additivity (\cite{K-T}). 
Given a sequence of $\O_X$ modules $A\subset B\subset C$ such that the multiplicity of the pairs is well defined, 
then $$e(A,C)=e(A,B)+e(B,C).$$
The result follows by setting $A=JM(X,f)$, $B=JM(X)\oplus \O_X$, $C=H_{d-1}(X)\oplus \O_X$, and using the fact that
$\?{JM(X,L)}=JM(X)\oplus \O_X$, so the multiplicity of the pairs $JM(X,L), H_{d-1}(X)\oplus \O_X$ and 
$JM(X))\oplus \O_X, H_{d-1}(X)\oplus \O_X$ are the same.

As the proposition shows, if we use $e(JM(X,f), JM(X)\oplus \O_X,x)$, or 

\noindent$(JM(X,f), H_{d-1}(X,0)\oplus \O_X,x)-(JM(X,L)),
H_{d-1}(X)\oplus \O_X,x)$ for our invariant, the number is the same. The second invariant shows that our number is measuring the extent to which the singularity of $f$ is more complicated than the singularity of a generic linear projection. 

Specialializing to the ICIS case shows an advantage of working with the difference of multiplicities, as these are easier to relate to other invariants of $X$; these relations give some insight into the meaning of our number. In the proposition below $\mu$ denotes Milnor number.

\prp 6 Suppose $X^d,0$ is an ICIS, $f$ a function with an isolated singularity, $L$ a linear function with no singularity at $0$.
Then $$e(JM(X,f), H_{d-1}(X)\oplus \O_{X},0)=\mu(X)+\mu(f),$$
 while 
$$e(JM(X,L), H_{d-1}(X)\oplus \O_X,0)=\mu(X)+\mu(X\cap L^{-1}(0)).$$

Hence
$$e(JM(X,f),JM(X)\oplus\O_X,0)=$$
$$e(JM(X,f), H_{d-1}(X,0)\oplus \O_X,0)-e(JM(X,L), H_{d-1}(X)\oplus \O_X,0)$$
$$=\mu(f)-\mu(L).$$

\pf Off the origin, $JM(X)$ has maximal rank, thus $H_{d-1}(X)=\O^p_X$, since $JM(X)=\O^p_X$ off the origin. Then
$e(JM(X,f), \O^{p+1}_X,0)=e(JM(X,f))=$ colength $J(JM(X,f))$. The last equality is 
by theorem 1.2 as $JM(X,f)$ has $n=d+(p+1)-1$ generators. The colength of $J(JM(X,f))$ by the theorem of L\^e and Greuel (\cite G, \cite L) is just 
$\mu(X)+\mu(f)$. The same argument applied to $JM(X,L)$ produces the second equality.

Looking at $e(JM(X,f), H_{d-1}(X)\oplus \O_X,0)$  in light of the ICIS experience
we expect that one part of the number is due to the
singularities of $X$. They enter in the $H_{d-1}(X)$ term and in the contribution of $X$ to $JM(X,f)$. Since $L$ has no singularity, 
(or minimal singularity depending on your perspective) the contributions from $X$ should cancel in the difference,
 which is what happens in 
the ICIS case.

Now we give an example of our ideas in the case where the singularities of $X$ are non-isolated. Suppose $X^d,0$ is a hypersurface defined
by a function $g$, the singular locus of $X$, $S(X)$, consists of an ICIS defined by an ideal $I$. Suppose the restriction of 
$g$ to a transverse slice to $S(X)$ is a Morse singularity, suppose $f$ is a
function on $X$  with an isolated singularity. Then $H_{d-1}(J(g))=I$. For it is clear that $J(g)\subset I$, $I$ is integrally closed
because it is radical, and by hypothesis on the singularity type of $g$, $J(g)=I$ except perhaps at $0$. So, our invariants become $e(JM(X,f),I\oplus\O_X)$ and $e(JM(X,L),I\oplus\O_X)$, where $L$ is a generic linear form.

The next proposition shows that if $f$ has no singularity at $x$ then our invariants cannot tell the difference between $f$ and a generic
linear form. This is important because if we deform our setup so that at the singular point of $X$, the deformed $f$ has no singularities,
then we may assume it is linear after that.

\prp 7 If $f$ has no singularity at $x\in X$, then $$e(JM(X,f),H_{d-1}(X)\oplus \O_X)=e(JM(X,L),H_{d-1}(X)\oplus \O_X)$$ \noindent where
$L$ is any linear function which does not have a singular point at $x$.

\pf If $f$ has no singularity at $x\in X$, then $JM(X,f)$ is  a reduction of $JM(X)\oplus \O_X$ as is $JM(X,L)$. So
 $$e(JM(X,f),H_{d-1}(X)\oplus \O_X)=e(JM(X)\oplus \O_X,H_{d-1}(X)\oplus \O_X)$$
$$= e(JM(X,L),H_{d-1}(X)\oplus \O_X).$$

Note that it may be the case that $\?{JM(X)}=H_{d-1}(X)$ at $x$. This means that if $L$ is a generic linear form for $X$, then 
 $e(JM(X,L),H_{d-1}(X)\oplus \O_X,x)=0$. For $\?{JM(X,L)}$ contains $JM(X)_L\oplus \{0\}$, hence  $\?{JM(X)_L}\oplus \{0\}$, hence
$H_{d-1}(X)\oplus \{0\}$, from which it follows that $\?{JM(X,L)}=H_{d-1}(X)\oplus \O_X$.

The next step is to construct a deformation of $f$ which will split  $df\cap T^*(X)$ into tranvserse intersections points. Our lemma is
similar to that of \cite{Le1}. The difference is in the notion of isolated singularity. Roughly, in \cite{Le1}, it is assumed that the
restriction of $f$ to each stratum has an isolated singularity.

\lem 8  Suppose $X,0\subset\IC^N$ is the germ of an equidimensional complex analytic set, $f$ the germ of a function on $X,0$ with an
isolated singularity at $0$. Then there is an open dense set $\Omega$ in the set of linear forms of $\IC^n$ such that for any form
$L\in\Omega$, $f+tL$ has only Morse singularities close to $0$ for $t$ sufficiently small, and no singularity in $S(X)$.

\pf The proof is very similar to Theorem 2.2 of \cite{Le1}, so we just sketch it. Let $\Theta:T^*(X)\times(\IC^N)^*\to (\IC^N)^*\times
(\IC^N)^*$ be defined by 
$$\Theta((z,k),l)=(k-(df(z)+l),l).$$

 Then $\Theta^{-1}(0\times(\IC^N)^*)$ is isomorphic to $T^*(X)$. The condition
that $f$ has an isolated singularity at $0$ is equivalent to the projection of $\Theta^{-1}(0\times(\IC^N)^*)$  to $\IC^N)^*$ being finite.
Now the desired $L$ can be found by chosing an element of $\IC^N)^*$ which is not in the image of $T^*(X)|S(X)$ nor in the discriminant of
the projection.

Now we are ready to prove our main result of this section.

\thm 9 Suppose $f:X^d,0\to\IC,0$ has an isolated singularity at the orign. Then 
$$\Gamma (df)\cdot T^*(X)=e(JM(X,f),JM(X)\oplus \O_X,0)$$
$$=e(JM(X,f),H_{d-1}(X)\oplus \O_X,0)-e(JM(X,L),H_{d-1}(X)\oplus \O_X,0)$$

\noindent where $L$ is a linear function without singularity at the origin.

\pf  Consider the family of sets $X\times\IC$. Let $\pi_X$ denote the projection from $X\times\IC$ to $X$. Choose $L$ as in the previous
lemma, then we have a family of functions $f_t=f+tL$ and associated family of graphs $\Gamma(d(f+tL))$ in $\pi^*_X(\IC^{n*})$. By
conservation of number,
$\Gamma (df)\cdot T^*(X)$ is  just $\Gamma (df+tL)\cdot T^*(X)$ for $t$ close to $0$, and this is just the number of Morse points of
$f+tL$.

We will compute the number of Morse points using the multiplicity-polar theorem.  The family of sets is
$X\times\IC$, the family of functions as above.
The two families of  modules are $JM(X,f_t)_{\pi_t}$ and $H_{d-1}(X) \O_{X\times\IC}\oplus \O_{X\times\IC}$.

First we must show that $\?{JM(X,f_t)_{\pi_t}}$ and $H_{d-1}(X) \O_{X\times\IC}\oplus \O_{X\times\IC}$ agree off an analytic set finite
over $\IC$. We know that the projection of $\Gamma (df+tL)\cdot T^*(X)$, to $X\times\IC$ defines the germ of an analytic set  at $0$ which
is finite over $\IC$. We claim that off this set $\?{JM(X,f_t)_{\pi_t}}=H_{d-1}(X) \O_{X\times\IC}\oplus \O_{X\times\IC}$.

Note that  $\pi^*_X(T^*X)=T^*(X\times\IC)$, so $\Gamma (df+tL)$ misses 
$T^*(X)$ at $(z,t)$ is equivalent to saying that no limiting tangent hyperplane to $X\times\IC$ at $(z,t)$ contains the kernel in $\IC^n$ 
of $d(f+tL)(z)$. In turn, as in lemma 2.1, 2.2 and 2.3, this implies the claim.

 We must also show that 
 $$C(\Projan {\cal R}(JM(X,f_t)_{\pi_t}))(0)=C(\Projan {\cal R}(JM(X,f))),$$ 
\noindent except possibly over $(0,0)\in X\times 0$.

Since $N=H_{d-1}(X) \O_{X\times\IC}\oplus \O_{X\times\IC}$ as a family of modules is independent of $t$, $\Projan {\cal R}(N)$ is a product, hence
$C(\Projan {\cal R}(N))(0)=C(\Projan {\cal R}(N(0)))$. Now, by the claim, at any point $p$ of $X\times 0$ close to the origin, there exists
a neighborhood $U$ of $p$ such that on $U$, $\?{JM(X,f_t)_{\pi_t}}=N$. This implies that over $U$, $\Projan {\cal R}(N)$ is finite over 
$\Projan {\cal R}(JM(X,f_t)_{\pi_t})$, and on $U\cap X\times 0$, $\Projan {\cal R}(N(0))$ is finite over 
$\Projan {\cal R}(JM(X,f_0))$. Now, since $\Projan {\cal R}(JM(X,f_0))\subset \Projan {\cal R}(JM(X,f_t)_{\pi_t})(0)$, the desired equality follows,
for any element of  $\Projan {\cal R}(JM(X,f_t)_{\pi_t})(0)$ has a preimage in  $\Projan {\cal R}(N))(0)$ which is  $\Projan {\cal R}(N(0))$, and the last
set maps to
$\Projan {\cal R}(JM(X,f_0))$. So, the multiplicity polar theorem applies. 

Note, that since $\Projan {\cal R}(N)$ is a product, $N$ has no polar curve. Now $JM(X,f_t)_{\pi_t}$ has $n$ generators, while the dimension
of $\Projan {\cal R}(JM(X,f_t)_{\pi_t})$ is the dimension of the base plus the generic rank of $D(F,f_t)$ minus $1$, where $F$ defines $X$. Thus
we have 

$${\rm{dim}} \Projan {\cal R}(JM(X,f_t)_{\pi_t})=d+1+(n-d)=n+1$$

Since the dimension of $\Projan {\cal R}(JM(X,f_t)_{\pi_t})$ is greater than or equal to the number of generators, there is no polar curve for
$JM(X,f_t)_{\pi_t})$.

It remains to determine the support of $ N/\?{JM(X,f_t)_{\pi_t}}$. At the Morse points of $f+tL$, $X$ is smooth, so $N$ is isomorphic to 
$\O^{n-d+1}_{X}$, and these points are in the support as $D(F,f_t)$ has less than maximal rank. At singular points of $X$, if $x\ne 0$, then
for $x$ close to $0$, $N(t)=\?{JM(X)}\oplus\O_{X,x}$. Since $f+tL$ has no singularity at $x$, it follows that $ N/\?{JM(X,f_t)_{\pi_t}}=0$.
So, the only point of $S(X)$ in the cosupport is $0$.

Then, the multiplicity polar theorem implies

$$e(JM(X,f),H_{d-1}(X)\oplus \O_X,0)=$$
$$\Gamma (df)\cdot T^*(X)+e(JM(X,f_t),H_{d-1}(X)\oplus \O_X,0).$$

Finally, since $f_t$ has no singular point at $0$ for $t$ small, we can replace it by a generic linear form, by proposition 2.7. So we get:

$$e(JM(X,f),H_{d-1}(X)\oplus \O_X,0)=\Gamma (df)\cdot T^*(X)+e(JM(X,L),H_{d-1}(X)\oplus \O_X,0).$$

From which the theorem follows.

\cor 10 The defect in the formula of \cite {BMPS} in the case of isolated singularities is

$$ (-1)^d (e(JM(X^d,f),H_{d-1}(X^d)\oplus \O_X^d,0)-e(JM(X^d,L),H_{d-1}(X^d)\oplus \O_X^d,0))$$
\pf By \cite{BMPS}, the defect is $(-1)^d \Gamma (df)\cdot T^*(X)$.

\sctn {Further applications and examples}

As was stated in the introduction, the number $\Gamma (df)\cdot T^*(X)$ has an importance that goes beyond the defect of the Euler
obstruction. In this section we look at some examples of this.

In a series of papers,(\cite {E-GZ1},\cite{E-GZ2},\cite{E-GZ3}, and \cite{EGZS})
 Ebeling and Gussein-Zade have been looking at the index of a differential 1-form $\omega$ defined on a singular space. In the case that
that $X$ has an isolated singularity they defined the radial index of a differential 1-form $\omega$ at a singular point $z$
denoted ${\rm{ind}}_z(\omega)$, in \cite{E-GZ1}.
The radial index has the property that on smooth varieities it agrees with the usual notion of the index, and satisfies conservation of
number.  Applying the result of the work of the last section, the natural module to use is $JM(X,\omega)$, which is the module
generated by the columns of the matrix whose last row is $\omega$, and whose first rows come from $DF$. Let $\tilde\chi(M)$
 denote the
reduced Euler characteristic of $M$. Based on the results  of \cite{EGZS} we can
show: 

\thm 1 Suppose $X^d,0$ is a germ of complex analytic set with an isolated singularity at $0$. Suppose $\omega$ is a differential 1-form with
an isolated singularity at $0$. Suppose $L$ is generic linear form with respect to the singularity of $X$. Then

$${\rm{ind}}_0(\omega)=e(JM(X,\omega), JM(X)\oplus \O_X,0)+(-1)^{d-1}\tilde \chi(L^{-1}(t)\cap X)$$

$$=e(JM(X,\omega), H_{d-1}(X)\oplus \O_X,0)-e(JM(X,dL), H_{d-1}(X)\oplus \O_X,0)$$
$$+(-1)^{d-1}\tilde \chi(L^{-1}(t)\cap X).$$

\pf From conservation of number and the existence of a generic $L$ which is dealt with in \cite{EGZS}, it follows that

$${\rm{ind}}_0(\omega)=C+{\rm{ind}}_0(dL),$$
\noindent where $C$ is the number of simple singularities of $\omega+tdl$ for small $t$. In turn this number is just
$\Gamma (\omega)\cdot T^*X$, where $\Gamma (\omega)$ denotes the image of the section of $T^*\IC^n|X$ induced by $\omega$.
By the proof of Theorem 2.9, this is just $e(JM(X,\omega), H_{d-1}(X)\oplus \O_X,0)-e(JM(X,dL), H_{d-1}(X)\oplus \O_X,0)$. To finish the
proof, in \cite{EGZS}, it is shown that ${\rm{ind}}_0(dL)=(-1)^{d-1}\tilde \chi(L^{-1}(t)\cap X).$

Our next application is to the theory of D modules. In \cite {Le2}, L\^e studies holonomic D modules $M$ on a neighborhood $U$ of the origin in $\IC^n$ and
functions
$f:U\to\IC$. Associated to $M$ there is its characteristic variety, which consists of the conormals of a collection of varieties $\{X_{\alpha}\}$ in
$U$. There is also the characteristic cycle, denoted $\tilde V$ which is the cycle whose underlying set is the characteristic variety  where
the multiplicity of $X_{\alpha}$, denoted $m_{\alpha}$ is the complex dimension
of the space of microlocal solutions to the complex analytic linear differential system associated to $M$ at a general point of $X_{\alpha}$. Assuming that the intersection
of $\Gamma (df)$ and $\tilde V$ is isolated at the origin, L\^e shows that the intersection number $\tilde V\cdot \Gamma(df)$ is the complex
dimension of the space of vanishing cycles of $M$ relative to $f$ at $0$. Now we can compute the intersection number $\Gamma (df)\cdot\tilde
V$, using the multiplicity of pairs of modules.

\thm 2 Suppose $M$ is a complex analytic linear differential system on a neighborhood $U$ of the origin in $\IC^n$ and
$f:U\to\IC$. Suppose $\{X_{\alpha}\}$ is the collection of subvarieties of $U$ whose conormals are the underlying sets of the
characteristic cycle $\tilde V$ of $M$. Suppose the intersection of $\Gamma(df)$ and $\tilde
V$ is isolated, $L$ generic for $f$ and all $X_{\alpha}$. Then
$$\Gamma (df)\cdot\tilde V=\Sum_{\alpha}m_{\alpha}e(JM(X_{\alpha},f),JM(X_{\alpha})\oplus\O_{X_{\alpha},0})=$$
$$\Sum_{\alpha}m_{\alpha}(e(JM(X_{\alpha},f),H_{{d_{\alpha}}-1}(X)\oplus\O_{X_{\alpha},0})-e(JM(X_{\alpha},L),H_{{d_{\alpha}}-1}(X)\oplus\O_{X_{\alpha},0})).$$

\pf The proof of this result follows immediately from the proof of Theorem 2.9, since the underlying sets of the components of $\tilde V$ are the conormals
of the $X_{\alpha}$.

In the important paper \cite{Ma1}, Massey calculated the homology of the Milnor fiber of a function $f$ with a stratified isolated
critical point on an analytic set $X$. (The function $f$ has a stratified isolated critical point if the union of the critical points of the restriction
of
$f$ to each stratum of the canonical Whitney stratification of $X$ is an isolated point.) We can restate his result in terms of the
multiplicity of the pair. In what follows $F_{f,p}$ denotes the Milnor fiber of $f$ at $p$. $L_{\alpha}$ is the complex link of the 
stratum $S_{\alpha}$ in $X$, $cone(L_{\alpha})$ is the cone over the complex link. If $G$ is an abelian group, then the notation 
$G^k$ denotes the direct sum of $k$ copies of $G$.

\thm 3 Let $\{S_{\alpha}\}$ be a Whitney stratification of $X$, suppose that $f:X,0\to \IC,0$ has a stratified isolated critical point at
$0$, $L$ a linear form which is  generic for $f$ and all $S_{\alpha}$. Then

$$H^i(B_{\epsilon}\cap X,F_{f,p},Z)=\oplus_{\alpha}H^{i-d_{\alpha}}(cone(L_{\alpha}),L_{\alpha})^{k_{\alpha}}$$

\noindent where $d_{\alpha}=$ dim $S_{\alpha}$, $(B_{\epsilon}$ a suitably chosen ball about $0$, and 
$k_{\alpha}=1$ if $d_{\alpha}=0$ and otherwise
$$k_{\alpha}=e(JM(S_{\alpha},f), JM(S_{\alpha})\oplus\O_{S_{\alpha},0})$$
$$=(e(JM(S_{\alpha},f),H_{{d_{\alpha}}-1}(S_{\alpha})\oplus\O_{S_{\alpha},0})-e(JM(S_{\alpha},L),H_{{d_{\alpha}}-1}(S_{\alpha})\oplus\O_{S_{\alpha},0})).$$

\pf This result is a modification of Theorem 3.2 of Massey (\cite{Ma1}). In his paper Massey calculates the $k_{\alpha}$
using the relative polar curve of $S_{\alpha}$ with respect to $f$ and the generic linear form $L$. He points out that $k_{\alpha}$ is
nothing other than the number of critical points on $S_{\alpha}$ of the pertubation of $f$ by $L$. From this observation and the proof of 
Theorem 2.9, the result follows.

Massey's result shows how the information from the complex links, which is information about the geometry of $X$ expressed in terms 
of a generic linear form, should be supplemented to take into account the difference between $f$ and a generic linear form. (The terms with $d_{\alpha}\ne 0$ are the supplementary terms.) Since the terms with $L$ in them reflect the geometry of $X$, it  would be an
improvement in this description to replace the terms with $L$, with terms more explicitly linked with the geometry of $X$. We next turn to an example which shows what is possible.

We return to the setup of the example before proposition 2.7. Given $X,0$ a hypersurface in $\IC^n$, $f:X,0\to\IC,0$, we know that the
Milnor fiber of $f$ is a bouquet of $n-2$ spheres. We assume that $S(X)$ is an ICIS defined by an ideal $I$, and that the transverse
singularity type of $X$ is Morse at all points except the origin.  We will study this situation by constructing a good deformation of it, then
using the multiplicity-polar theorem to control the deformation. This will require some basic data from the case where $X$ has a $D_{\infty}$
singularity.  

\lem 4 Suppose $X^{n-1},0\subset\IC^n,0$ has a $D_{\infty}$ singularity at the origin, $X$ defined by $g$. Then
$$e(JM(X,l),I\oplus \O_X,0)=e(J(g)_l,I,0)=1.$$

\pf The assumption implies that the singular set of $X$ is a smooth curve and that $g$ is the suspension, up to a change of coordinates, by adding
disjoint square terms of $z^2-x^2y$. We'll assume $g(x,y,z)=z^2-x^2y$, as the extra square terms won't significantly affect the proof.

First we'll show that $e(J(g)_l,I,0)=1$, then see how the second computation reduces to the first.

 We have to compute a sum of 
intersection 
numbers:

$$e(J(g)_l,I,0)=\Sum_{j=0}^{1} \int D_{J(g)_l,I} \cdot l^{1-j}_{J(g)_l}\cdot 
l^j_{I}.$$

Consider the term $ \int D_{J(g)_l,I} \cdot
l_{I}$. The intersection $B\cdot l_{I}$ is a curve; this curve projects to $\Projan {\cal R}(I)$ and to $X$. If we want to calculate its image 
on $X$, we take the polar curve of $I$ on $X$. This is empty, because $V(I)$ is a set theoretic hypersurface. Hence $ \int D_{J(g)_l,I} \cdot
l_{I}=0$.

Now consider $ \int D_{J(g)_l,I} \cdot l_{J(g)_l}$. The inclusion of $J(g)_l$ in $I$ induces an inclusion of Rees algebras; in turn this induces an
ideal sheaf on
$B_I(X)$. For the example on hand, we can take $l$ to be a linear function whose kernel is not limiting tangent hyperplane plane to
$X$ at $0$ so $l=x$. Then $J(g)_l=<x^2,z>$. The Rees algebra of $I$ is isomorphic to $\O_X[T_1,T_2]/R_I$ where $R_I$ is the kernel of the map
from $\O_X[T_1,T_2]$ to ${\cal R}(I)$, where $T_1$ goes to $x$ and $T_2$ goes to $z$. Then the ideal sheaf induced by the inclusion of Rees algebras is 
$(xT_1/T_2, 1)$ when $T_2\ne 0$ and $(x, T_2/T_1)$ when $T_1\ne 0$. Now, as $(s,t^2,st)$ is a parameterization of $X$, $(s,t^2,st, [1,t])$ is a
parameterization of $B_I(X)$. Our ideal sheaf pulls back to $(s,t)$; calculating  $ \int D_{J(g)_l,I} \cdot l_{J(g)_l}$ here gives $1$, the multiplicity
of
$(s,t)$.

Now consider $e(JM(X,l),I\oplus \O_X,0)$. Here $\Projan {\cal R}(I\oplus \O_X)$ is obtained from $B_I(X)\subset X\times \IP^1$ by taking the join with a
fixed point
$[0,0,1]$ in $\IP^2$ of each of the fibers of $B_I(X)$ over $X$. There are three intersection numbers to compute; two of which involve a factor of
$l_{I\oplus\O_X}$. The intersection of $\Projan {\cal R}(I\oplus \O_X)$ with a generic hyperplane is isomorphic to $B_I(X)$; this time however, the ideal
sheaf is the constant sheaf. It remains to consider $ \int D_{JM(X,l)} \cdot l^2_{JM(X,l}$. Now the Rees algebra $I\oplus \O_X$
is isomorphic to $\O_X[T_1,T_2,T_3]/R_I$, where $T_3$ goes to $(0,1)$; the ideal sheaf induced from $JM(X,l)$ in the chart where $T_1\ne 0$ is 
$(y+T_3/T_1, x, T_2/T_1)$ (In the other charts the sheaf is the constant sheaf.)

We again have a parametrization of $\Projan {\cal R}(I\oplus \O_X)$ given by 
$$F: \IC^2\times\IP^1\to \Projan{\cal R}(I\oplus \O_X)$$
$$F(s,t,[S_1,S_2]=(s,t^2,st,[S_1,S_1t,S_2].$$ In the chart where $S_1\ne 0$, the ideal sheaf induced from $JM(X,l)$ pulls back to
$(s,t+S_2/S_1,S_2/S_1)$, so again the intersection number is 1. In the other chart, the sheaf pulls back to the constant sheaf. This finishes the
computation.

Note that inducing the  ideal sheaf from the inclusion of Rees algebras has the effect of ``dividing" $JM(X,l)$ by $I\oplus \O_X$.

The next proposition proves a useful relation between some of our invariants. In this proposition $\lambda^0(g)$ denotes the L\^e number of dimension 0.(Cf.\cite{Ma} for information about the L\^e numbers.)

\prp 5 In the above setup if the dimension of $S(x)$ is $1$ then $$e(JM(X,l),I\oplus \O_X,0)$$
$$=e(J(g)_l,I,0)=\lambda^0(g)+\mu(X\cap
l^{-1}(t))-2D_{\infty}(g)-2e(JM(V(I)),0)$$

Otherwise
$$e(JM(X,l),I\oplus \O_X,0)=e(J(g)_l,I)=\lambda^0(g)+\mu(X\cap H_t)-2e(JM(V(I)),0)$$

\pf  If the dimension of $S(X)$ is $1$, then we can construct a deformation of $g$, and hence of $X$, $I$, over a smooth base $Y$, so that for the generic
parameter value $S(X)$ is smooth, $X$ has only  of $D_{\infty}$ singularities and there $D_{\infty}(g)$ of them, and the only singularities of $l$ on
$X_t-S(X_t)$ are Morse. We can do more. We can ensure that $l$ has only Morse singularities on $S(X)$, and that $l^{-1}(0)$ is not a limiting tangent
hyperplane to $X_t$ at such points. Now the support of $(I_t\oplus \O_{X_t})/JM(X(t),L)$ and $I(t)/J(g_t)_l$ for generic $t$ values are the same. They
will consist of the $D_{\infty}$ points and the singularities of of $l$ restricted to $X(t)-S(X(t))$. By the multiplicity-polar theorem since neither 
$JM(X,l)$,$I\oplus \O_X$, $J(g)_l$ nor $I$, have a polar variety the same dimension as $Y$, and the multiplicities of the pairs have the same values
at the $D_{\infty}$ points it follows that the values of the multiplicities of the pairs are the same. If the dimension is greater than 1, then a good
deformation still exists but there are no $D_{\infty}$ points. (If there were there would be a hypersurface of points in $S(X)$ where the transversal
singularity type was not Morse.) This proves the first part of the two equations.

To prove the second part, we consider the intersection number $\Gamma^1(g)\cdot X$, where we use $l$ to define the relative polar curve, and how this
number breaks up under deformations. By Massey, the number is $\lambda^0(g)+\mu(X\cap l^{-1}(t))$. On the other hand, under deformations the intersection
product breaks into three different types of points: the singular points of $l$ on
$X_t-S(X_t)$ which count once, the singular points of $l$ restricted to the generic fiber of the miniversal deformation of $V(I)$, which count twice and of
which there are $e(JM(V(I))$, and the $D_{\infty}$ points which count three times.

These considerations give the equations:

If the dimension of $V(I)$ is 1, then

$$\lambda^0(g)+\mu(X\cap
l^{-1}(t))= 3D_{\infty}(g)+2e(JM(V(I)),0)+(e(J(g)_l,I,0)-D_{\infty}(g)$$

If the dimension of $V(I)>1$, then

$$\lambda^0(g)+\mu(X\cap
l^{-1}(t))= 2e(JM(V(I)),0)+(e(J(g)_l,I,0))$$

Solving for $e(J(g)_l,I,0)$ gives the desired results.

Now we are ready to prove our improvement on Massey's result in the current set up.

\thm 6 Suppose $X^{n-1}$ is a hypersurface in $\IC^n$ defined by $g$, and $S(X)$ is an ICIS defined by $I$, $X$ has transverse $A_1$ type along $V(I)$
except perhaps at the origin. Suppose $f:X,0\to \IC,0$ has a stratified isolated critical point. Then Milnor fiber of $f$ has the homotopy type of a
bouquet  of
$n-2$ spheres, where the number of spheres in the bouquet, $b_{n-2}(f)$ is $$e(JM(X,f),I\oplus
\O_X,0)-(e(J(g),I,0;\O_n)-D_{\infty}(g))+e(JM(V(I),f,0).$$

\noindent if the dimension of $S(X)=1$, and otherwise

$$b_{n-2}(f)=e(\?{JM(X,f)},I\oplus
\O_X,0)-e(J(g),I,0;\O_n)+e(JM(V(I),f,0).$$
\pf By a theorem of L\^e, we know that the  Milnor fiber of $f$ has the homotopy type of a
bouquet  of
$n-2$ spheres. We can get the number of spheres from Massey's result. Here the canonical stratification of $X$ has three strata and is $\{X-S(X), S(X)-0,
0\}$.

Assuming
$f$ has a stratified isolated singularity, Massey's result gives that the number of spheres $b_{n-2}(f)$ is 

$$b_{n-2}(f)=(e(JM(X,f),I\oplus\O_X)-e(JM(X,l),I\oplus\O_X))$$
$$+(e(JM(V(I),f)-e(JM(V(I),l))+\mu(X\cap l^{-1}(t))$$
 
Here the term in the first parenthesis comes from the open stratum, the second parenthesis from $V(I)-0$, and the third term from $0$.
Now we use proposition 3.5, and we get if the dimension of $V(I)$ is 1:

$$b_{n-2}(f)=(e(JM(X,f),I\oplus\O_X)$$
$$-(\lambda^0(g)+\mu(X\cap
l^{-1}(t))-2D_{\infty}(g)-2e(JM(V(I)),0)$$
$$+(e(JM(V(I),f)-e(JM(V(I),l))
+\mu(X\cap l^{-1}(t)).$$

From \cite {G-7} we have a formula for $\lambda^0(g)$, if the dimension of $V(I)$ is 1,

$$\lambda^0(g)=e(J(g),I,0;\O_n)+e(JM(V(I))+D_{\infty}(g).$$

The same formula holds without the $D_{\infty}(g)$ term if the dimension of $V(I)$ greater than 1.

Substituting for $\lambda^0(g)$ and canceling gives: 

$$b_{n-2}(f)=(e(JM(X,f),I\oplus\O_X)-(e(J(g),I,0;\O_n)-D_{\infty}(g))+e(JM(V(I),f)).$$
In the case where dimension of $S(X)>1$, there is no $D_{\infty}(g)$ term and we get the formula for this case.

Notice that the dependence on our linear form has disappeared. If $X$ has an isolated singularity, then the only terms that don't drop out are
$e(JM(X,f),I\oplus\O_X)$ and $e(J(g),I,0;\O_n)$. These become $e(JM(X,f),\O^2_X)$ and $e(J(g),0;\O_n)$, since $I=\O_X$ or $\O_n$. Since
$e(JM(X,f),\O^2_X)=\mu(f)+\mu(X)$, and we have that $e(J(g),0;\O_n)=\mu(X)$, the difference is $\mu(f)$ as it should be.

We have chosen to group together the $e(J(g),I,0;\O_n)$, $D_{\infty}(g)$ terms in our
formula because their difference is the number of Morse points coming out of a deformation of $g$.

 We mention one more example, which gives some insight into the meaning of our result. If $g=z^2-x^2y$ and $f=y-z$, then $f$ is a non-generic linear function
which  has a stratified isolated singularity. In this case, the $(e(J(g),I,0;\O_n)-D_{\infty}(g))$ term is zero, the $e(JM(V(I),f)$ is zero and the
$e(JM(X,f),I\oplus\O_X)$ term is 2 by the same method of calculation as for a generic linear form in the $D_{\infty}$ case. So the Milnor number of $f$ on
$X$ is 2. This can be checked as follows: The two equations $g=f=0$ defines an isolated singularity $z(z-x^2)=0$. This has Milnor number 3--however the
Milnor fiber here is defined by the intersection of $g=\epsilon_1$ with $f=\epsilon_2$. Meanwhile, the curve defined by the intersection of $g=0$ with
$f=\epsilon_2$ has a Morse singularity at the origin. Thus deforming from the Milnor fiber of $f|g=0$ to the Milnor fiber of $(f,g)$ adds 1 vanishing cycle
for a total of three, so the Milnor fiber of $f|g=0$ has 2 vanishing cycles as predicted. Since $X$ has non-isolated singularities, the Milnor fiber of $f$
on $X$ has singularities.

Our last application is to the \Af stratification condition. In a previous paper \cite{G-6}, we described how to get a necessary and sufficient criterion in
terms of analytic invariants for the \Af condition.
The viewpoint of this paper allows us to reformulate this  criterion in a way which is less dependent
on the family of sets.  We will be working with families of sets with isolated singularities. If
$X^d,0$ is analytic set germ with an isolated singularity, then $H_0(JM(X),0)= H_{d-1}(JM(X),0)$. For $h\in H_0(JM(X),0)$ implies $h\in JM(X)$ off a set of
codimension $1$, since $JM(X)$ is integrally closed at smooth points. But it is easy to see that if $h\in JM(X)$  off a set of codimension $1$, then it is in 
$JM(X)$ off the origin, since $X$ has an isolated singularity.

There is a geometric description of 
$H_0(JM(X),0)$, if $X$ is not a complete intersection. We have $JM(X)\subset \O^p_X$ for some $p$, hence  $JM(X)$, induces an ideal sheaf ${\cal M}$, on $\Projan {\cal R}(\O^p_X)$. Since $X$ is not a complete intersection, there must be a component of $V({\cal M}))$ for each component of $X$ which surjects onto $X$ under projection to the base. Denote the union of these components by $V_0({\cal M})$.
Let $I_0(X)=I(V_0({\cal M}))$. Then the terms of degree 1 of $I_0(X)$ are the elements of $H_0(JM(X),0)$. For the points of $V_0({\cal M})$ over the smooth points of $X$ are hyperplanes which
contain the column space of a matrix of generators of $JM(X)$ at a point. The condition that $h\in JM(X)$ at a smooth point of $X$ is equivalent to asking
that $h$ lies in every hyperplane that contains the column space of a matrix of generators of $JM(X)$. Thus $h\in JM(X)$ at all smooth points of $X$ iff the induced element of ${\cal R}(\O^p_X)$ vanish on all points of $(V_0({\cal M}))$ over the smooth points of $X$ iff the induced element lies in $I_0(X)$.

There is also an alternate algebraic description of $H_0(JM(X),0)$. Denote the set of all elements  $h\in I$  such that the
partial derivatives of $h$ are in $I$ by
$\int I$. Note that we can identify $\O^p_X$ with its dual $\hom(\O^p_X, \O_X)$. If $I$ has $p$ generators, we have the following short exact sequence of
$\O_X$ modules.

$$0\To R\to \O^p_X\To I/I^2\To 0$$ 

Here $R$ is the module of relations. Denote the map to $I/I^2$ by $j$. This gives the injection 
$$0\to \hom(I/I^2,\O_X)\to \hom(\O^p_X, \O_X).$$

So, we can identify elements in the image of this last inclusion with their preimages. Note that each partial derivative operator defines an element of 
$\hom(I/I^2,\O_X)$. Denote the submodule of $\hom(I/I^2,\O_X)$ generated by the elements defined by the partial derivative operators by $D$.

\prp 7 With the identification of $\O^p_X$ with $\hom(\O^p_X, \O_X)$, the module $JM(X)$ is the image of $D$ under the inclusion of $\hom(I/I^2,\O_X)$ in
$\hom(\O^p_X, \O_X)$, and $H_0(X)$ is the image of $\hom(I/\int I,\O_X)$. 

\pf Note that applying $\pd{}{z_i}$ to the generators of $I$ produces a column vector of a matrix of generators of $JM(X)$. This column vector is the image of 
$\pd{}{z_i}$ in $\hom(\O^p_X, \O_X)$. The module $R$ is a submodule of larger module $R'$ which is the intersection of the kernels of the elements of 
$JM(X)$ viewed as a submodule of $\hom(\O^p_X, \O_X)$. The elements of $H_0(X)$, viewed as a submodule of $\hom(\O^p_X, \O_X)$ are exactly the 
elements of $\hom(\O^p_X, \O_X)$ whose kernel contains $R'$. This is clearly necessary and  among the elements of $R'$, there are p-tuples of maximal
non-vanishing minors. For these to be in the kernel of $h$, implies that the generic rank of the matrix of generators of $JM(X)$ augmented by $h$ is the same
as the generic rank of the matrix of generators of $JM(X)$, which implies $h$ generically in $JM(X)$.

Now $R'$ in the kernel of $h$ implies that $h$ is in the image of $\hom(I/I^2,\O_X)$, and further that the preimage of $h$ vanishes on $\int I$, hence
$h$ comes from an element of $\hom(I/\int I,\O_X)$.

 The \Af result from \cite{G-6} involves the multiplicity of the polar variety of codimension $d$ of $H_0(JM(X,F))$, which is one of the modules in the
pair.  Given a family of pairs of modules $(M(y),N(y))$, the invariant which incorporates the polar variety information of 
the larger is $e_{\Gamma}(M,N,y)$, defined
as follows:

$$e_{\Gamma}(M,N,y):=\Sum_{x\in p^{-1}(y)}e(M,N,x)+{\rm mult}_Y\Gamma(N,x)$$

It may be necessary to re-choose the representative of the polar variety at different $x$.

 The following result follows from
the material of
\cite{G-6}.

\thm 8   Suppose $( X^{d+k},0) \subset  (\IC^{n+k},0)$, 

$X = G^{-1}(0)$,  $G:{\IC}^{n+k}\to{\IC}^p$, $Y$ a smooth subset of $X$,
coordinates  chosen so
that ${\IC}^k \times {0} = Y$, $X$ equidimensional with equidimensional fibers, $X$ reduced, all fibers generically reduced, 
and fibers reduced over a Z-open subset of $Y$. 
Suppose $F\:(X,Y)\to (\IC,0)$,  $F\in m^2_Y$, $Z=F^{-1}(0)$. Suppose $X_y$ and $Z_y$ are isolated
singularities, suppose $c(JM(G,F))$ is $Y$.

Then $e_{\Gamma}(JM(G_y;F_y),H_0(JM_z(G;F))(y))$ is independent of $y$, 
iff the pair of strata
$(X-Y,Y)$ satisfies Thom's \AF \hskip 2pt condition.

\pf This follows from Theorem 5.6 of \cite{G-6}.

We can improve this result by removing the dependence on the polar variety of $ H_0(JM_z(G;F))$ by using a generic linear function.

\thm 9 Suppose $( X^{d+k},0) \subset  (\IC^{n+k},0)$, 

$X = G^{-1}(0)$,  $G:{\IC}^{n+k}\to{\IC}^p$, $Y$ a smooth subset of $X$,
coordinates  chosen so
that ${\IC}^k \times {0} = Y$, $X$ equidimensional with equidimensional fibers, $X$ reduced, all fibers generically reduced, 
and fibers reduced over a Z-open subset of $Y$. 
Suppose $F\:(X,Y)\to (\IC,0)$,  $F\in m^2_Y$, $Z=F^{-1}(0)$. Suppose $X_y$ and $Z_y$ are isolated
singularities, suppose $c(JM(G,F))$ is $Y$. Suppose $(X-Y,Y)$ satisfies Thom's \AL condition with respect to $L$ a generic linear function 
on $\IC^n$. Then
 $$e(JM(G_0;F_0),H_0(JM_z(G;F))(0),0)-e(JM(G_0,L),H_0(JM_z(G;L))(0),0)$$
$$=e(JM(G_y;F_y),H_0(JM(G_y;F_y),(0))-e(JM(G_y,L),H_0(JM(G_y;L),0)$$ \noindent for all $y$ generic close to $0\in Y$ 
iff the pair of strata
$(X-Y,Y)$ satisfies Thom's \AF \hskip 2pt condition.

\pf  The hypothesis on the \AL condition implies by Theorem 3.7 that the invariant $e_{\Gamma}(JM(G_y,L,),H_0(JM_z(G;L))(y), y)$ is independent of $y$. Hence the invariant

$$e_{\Gamma}(JM(G_y,F_y),H_0(JM_z(G;F))(y), y)- e_{\Gamma}(JM(G_y,L,),H_0(JM_z(G;L))(y), y)$$ 
\noindent is independent of $y$ iff 
$$e_{\Gamma}(JM(G_y;F_y),H_0(JM_z(G;F))(y))$$ 
\noindent is. Since $H_0(JM_z(G;F))
=H_0(JM_z(G;L))$, as both $F_y$ and $L$ have isolated singularities on each fiber,  the polar variety terms cancel in the difference, from which the result
follows.

Note in Theorem 3.9, the only remaining place where information from the fiber enters the invariant is in the $H_0(JM_z(G;F))(0)$ term which is part of the
invariant at the origin. The Z-open set referred to in the statement of the theorem is chosen so that for $y$ in this set
$H_0(JM_z(G;F)(y))=H_0(JM(G_y,F_y))$. For details on the existence of this Z-open set see lemma 5.2 of \cite{G-6}.

This observation allows us to link the discussion of the \AF condition back with the defect.

\cor 10 Suppose, in addition to the hypotheses of Theorem 3.9, we have $H_0(JM_z(G))(0)=H_0(JM(G_0))$. Then the \AF condition holds if and 
only if the defect of $f_y$ at the origin is independent of $y$.

\pf The additional hypothesis implies that $$H_0(JM_z(G;L))(y)=H_0(JM_z(G;F))(y)=H_{d-1}(JM(G_y))\oplus \O_X.$$ 
These equalities follow by an argument similar to that for lemmas 2.2 and 2.3 and using the fact that for isolated singularities $H_{d-1}(JM(G))=H_0(JM(G))$.

These equalities imply the difference of the multiplicities in Theorem 9 is just the defect up to a sign.

It is also possible to use $\hom(I/I^2,\O_X)$ to develop the theory instead of $H_0(X)$.  This has the advantage that if the base of the miniversal
deformation of $X$ is smooth, then it is obvious that $\hom(\I/\I^2,\O_{\cal X})(0)=\hom(I/I^2,\O_X)$. It is an interesting question to see what conditions
a deformation of $X$ must satisfy for this equality to hold in general.

\sct References

\references

 Br  
 J.P. Brasselet,{\it Existence des classes de Chern en th\'eorie bivariante}, Ast\'erisque,
vol. 101-102,  7--22, 1983

 BLS 
J.P. Brasselet,  L\^e D. T., J. Seade, { Euler obstruction and indices of vector fields}, Topology 
  vol. 39 ,  1193-1208, 2000

BMPS
 J.-P. Brasselet, D. Massey, A. J. Parameswaran and J. Seade,{\it Euler Obstruction and Defects of Functions on Singular Spaces} preprint,
2003

B-R
 D. A. Buchsbaum and D. S. Rim,
  {\it A generalized Koszul complex. II. Depth and multiplicity,}
 \tams 111 1963 197--224

E-GZ1
 W. Ebeling, S. M. Gusein-Zade: On the index of a vector field at an 
isolated singularity. In: The Arnoldfest, edited by E. Bierstone et al., 
Fields Inst. Commun. 24, AMS, 1999, pp. 141-152

E-GZ2
 W. Ebeling, S. M. Gusein-Zade, On the index of a holomorphic 1-form on an isolated complete 
intersection singularity. Doklady Math. 64 (2001), 221-224

E-GZ3
 W. Ebeling, S. M. Gusein-Zade: Indices of 1-forms on an isolated complete 
intersection singularity. Moscow Math. J. 3, 439-455 (2003)

EGZS
 W. Ebeling, S. M. Gusein-Zade, J. Seade: Homological index for 
1-forms and a Milnor number for isolated singularities. Preprint 
math.AG/0307239 

F
 W. Fulton,
 ``Intersection Theory,''
 Ergebnisse der Mathematik und ihrer Grenzgebiete, 3. Folge
 $\cdot$ Band 2, Springer--Verlag, Berlin, 1984

FM 
 W. Fulton and R. MacPherson, {\it Categorical Framework for the study of singular Spaces},
 Memoirs of Amer. Math. Soc.
vol. 243 1981

 G-1
   T. Gaffney,
   {\it Aureoles and integral closure of modules,}
   in ``Stratifications, Singularities and Differential
 Equations~II'',  Travaux en Cours {\bf 55}, Herman, Paris, 1997, 55--62

 G-2
   T. Gaffney,
   {\it Integral closure of modules and Whitney equisingularity,}
   \invent 107 1992 301--22

 G-3
   T. Gaffney,
   {\it Plane sections, \Wf and \Af,} in ``Real and complex singularities 
(S\~ao Carlos, 1998)'',
 Chapman and Hall Res. Notes Math. {\bf 412}, 2000, 17-32

 G-4
   T. Gaffney,
   {\it Multiplicities and equisingularity of ICIS germs,}
    \invent 123 1996 209--220

G-5
	T. Gaffney, {\it Generalized Buchsbaum-Rim Multiplicities and a 
Theorem of Rees,} Communications in Algebra, vol 31 \#8 p3811-3828, 2003

G-6 
 T. Gaffney,{\it Polar methods, invariants of pairs of modules and 
equisingularity,}  Real and Complex Singularities (Sao Carlos, 2002),
Ed. T.Gaffney and M.Ruas, Contemp. Math.,\#354, Amer. Math. Soc.,
Providence,
RI, June 2004, 113-136

G-7
T.Gaffney, {\it The multiplicity of pairs of modules and hypersurface 
singularities}, submitted to Real and Complex Singularities, (Sao Carlos-Luminy 2004)

G-8
	T. Gaffney, {\it The Multiplicity-Polar Formula and Equisingularity}, 
in preparation

GK
   T. Gaffney and S. Kleiman,
   {\it Specialization of integral dependence for modules}
   \invent 137 1999 541-574

G
 G. M. Greuel,
 ``Der Gauss--Manin Zusammenhang isolierter Singularit\"aten von
vollst\"and\-igen Durchschnitten," Dissertation, G\"ottingen (1973),
\ma 214 1975 235--66

 K-T
   S. Kleiman and A. Thorup,
    {\it A geometric theory of the Buchsbaum--Rim multiplicity,}
   \ja 167 1994 168--231

 KT1
 S. Kleiman and A. Thorup,
  {\it The exceptional fiber of a generalized conormal space,}
  in ``New Developments in Singularity Theory." D. Siersma, C.T.C. Wall and 
V. Zakalyukin (eds.), Nato Science series,
II Mathematics, Physics and Chemistry-Vol. 21 2001 401-404

L
 D. T. L\^e,
  {\it Calculation of Milnor number of isolated singularity of complete
intersection,}
 \faa 8 1974 127--31

Le1 
  D. T. L\^e, {\it Le concept de singularit\'e isol\'ee de function analytique}, Advanced studies in pure math. 8 (1986), 215-227

Le2
	 D. T. L\^e, {\it Morsification of  D modules}, Bol. Soc. Mat. Mexicana (3) 4 (1998), no. 2, 229--248.

 MP   
 R. MacPherson,{\it Chern classes for singular varieties}
 Ann. of Math, vol. 100, 423-432 1974

Ma
  D. Massey {\it L\^e Cycles and Hypersurface Singularities,}
  Springer Lecture Notes in Mathematics 1615 , (1995)

Ma1
	D.Massey {\it Hypercohomology of Milnor Fibers}, Topology, v 35, \#4,p969-1003,(1996)

Ma-2
  D. Massey {\it Numerical Control over Complex Analytic Singularities}
 Memoirs of the AMS,\#778, AMS 2003

 T-2
   B. Teissier,
   {\it Multiplicit\'es polaires, sections planes, et conditions de
 Whitney,}
   in ``Proc. La R\'abida, 1981.'' J. M. Aroca, R. Buchweitz, M. Giusti 
and
 M.  Merle (eds.), \splm 961 1982 314--491

\endreferences

\bye